\documentclass[final]{siamltex}

\usepackage{color}
\usepackage{a4,amsmath,amssymb,amscd,latexsym} 
\usepackage{bbm} 
\usepackage{epsfig,color} 
\usepackage{amsfonts}
\usepackage{mathrsfs}
\usepackage{overpic}
\usepackage{subcaption}

\usepackage{paralist}

\usepackage{graphicx}

\usepackage{trfsigns}
\usepackage{algpseudocode} 

\usepackage{url}

\usepackage{stmaryrd}

\newcommand{\Gi}{\mathrm{\Gamma_{int}}}
\newcommand{\Gb}{\mathrm{\Gamma_{bottom}}}
\newcommand{\Gl}{\mathrm{\Gamma_{left}}}
\newcommand{\Gr}{\mathrm{\Gamma_{right}}}
\newcommand{\Gt}{\mathrm{\Gamma_{top}}}
\newcommand{\Go}{\mathrm{\Gamma_{out}}}

\def\norm#1{\|#1\|} 

\def\nd#1{\frac{\partial #1}{\partial n}} 
\def\td#1{\frac{\partial #1}{\partial t}} 

\def\intdxdt#1{\int_{0}^{T} \int_{\Omega} #1 \, dx\, dt}

\def\intdx#1{\int_{\Omega} #1 \, dx}
\def\intdsidt#1{\int_{0}^{T}\int_{\Gi} #1 \, ds\, dt}

\def\intdsodt#1{\int_{0}^{T} \int_{\Go} #1 \, ds\, dt}

\def\intdt#1{\int_{0}^{T} #1 \, dt}

\def\grad{\mbox{grad}}
\newcommand{\LL}{{\mathscr{L}}}

\def\bei#1{\vrule width 0.4pt height 14pt depth 9pt
           \lower 8pt \hbox{$ _{\hbox{} #1}$}\!\!\!}

\newcommand{\ACS}{2pt}
\newcommand{\mat}[4]{{\arraycolsep\ACS
\left#1\begin{array}{@{}*{#2}{c}@{}}#4\end{array}\right#3}}

\newcommand{\N}{{\mathbbm{N}}} 
\newcommand{\R}{{\mathbb{R}}} 

\makeatletter
\def\moverlay{\mathpalette\mov@rlay}
\def\mov@rlay#1#2{\leavevmode\vtop{%
   \baselineskip\z@skip \lineskiplimit-\maxdimen
   \ialign{\hfil$\m@th#1##$\hfil\cr#2\crcr}}}
\newcommand{\charfusion}[3][\mathord]{
    #1{\ifx#1\mathop\vphantom{#2}\fi
        \mathpalette\mov@rlay{#2\cr#3}
      }
    \ifx#1\mathop\expandafter\displaylimits\fi}
\makeatother

\newcommand{\cupdot}{\charfusion[\mathbin]{\cup}{\cdot}}

\def\scp#1{\left\langle #1\right\rangle}
\def\norm#1{\left\| #1\right\|}

\def\CROP#1{}

\newtheorem{remark}{Remark}

\newenvironment{rev}  {\color{black}}{\color{black}}

\title{Efficient PDE constrained shape optimization based on Steklov-Poincar\'e type metrics}

\author{Volker H.~Schulz, Martin Siebenborn and Kathrin Welker\thanks{University of Trier, Department of Mathematics, 54296 Trier, Germany ({\tt volker.schulz@uni-trier.de, siebenborn@uni-trier.de, welker@uni-trier.de}).}
        }

\begin{document}
\begin{rev}
\maketitle
\end{rev}
\begin{abstract}
Recent progress in PDE constrained optimization on shape manifolds is based on the Hadamard form of shape derivatives, i.e., in the form of integrals at the boundary of the shape under investigation, as well as on intrinsic shape metrics. From a numerical point of view, domain integral forms of shape derivatives seem promising, which rather require an outer metric on the domain surrounding the shape boundary. This paper tries to harmonize both points of view by employing a Steklov-Poincar\'e type intrinsic metric, which is derived from an outer metric. Based on this metric, efficient shape optimization algorithms are proposed, which also reduce the analytical labor, so far involved in the derivation of shape derivatives. 
\end{abstract}



\pagestyle{myheadings}
\thispagestyle{plain}

%
\section{Introduction}

Shape optimization is of interest in many fields of application -- in particular in the context of partial differential equations (PDE). As examples, we mention aerodynamic shape optimization \cite{AIAA-2013}, acoustic shape optimization \cite{Berggren-horn-2007} or optimization of interfaces in transmission problems \cite{Langer-2015,Skin-2015,Paganini} and in electrostatics \cite{Zapletal_2015}. In industry, shapes are often represented within a finite dimensional design space. However, often this reduction is felt as being too restrictive \cite{Bletzinger-2015}, which motivates shape optimization based on shape calculus. Major effort in shape calculus \cite{Delfour-Zolesio-2001,SokoZol} has been devoted towards expressions for shape derivatives in so-called Hadamard-form, i.e., in boundary integral form. It is known that the second order shape derivative, formerly coined as shape Hessian, is nonsymmetric in general, which for a long time has been an obstacle for algorithmic developments in shape optimization in the fashion of nonlinear programming. Recently \cite{VHS-shape-Riemann,Schulz-Structure-2014,Schulz-LN-2014}, shape optimization has been considered as optimization on Riemannian shape manifolds, which enables design and analysis of NLP-like algorithms including one-shot sequential quadratic programming and theoretical insights into the structure of the second order shape derivative in comparison to the Riemannian shape Hessian. 
\begin{rev}
Coercivity results for shape Hessians for elliptic problems can be found in \cite{EpplerHarbrechtSchneider-2007}. The scalar product used in this work is in line with these results.
\end{rev}

On the other hand, it is often a very tedious, not to say painful, process to derive the boundary formulation of the shape derivative. Along the way, there frequently appears a domain formulation in the form of an integral over the whole domain as an intermediate step. Recently, it has been shown that this intermediate formulation has numerical advantages \cite{Berggren,Langer-2015,HipPag_2015,Paganini}. In \cite{LaurainSturm2013}, also practical advantages of the domain shape formulation have been demonstrated, since it requires less smoothness assumptions. Furthermore, the derivation as well as the implementation of the domain integral formulation requires less manual and programming work. Thus, there arises the natural goal of combining the favorable domain integral formulation of shape derivatives with the favorable NLP-type optimization strategies on shape manifolds, which seem so far tightly coupled with boundary integral formulations of shape derivatives. This publication aims at demonstrating that this coupling is indeed possible and that it naturally leads to a novel family of Poincar\'e-Steklov type metrics on shape manifolds.
\begin{rev}
In contrast to \cite{Schulz-Structure-2014} this work consciously avoids surface formulations of shape derivatives in order to provide more handy optimization algorithms.
\end{rev}

The paper is organized in the following way. 
First, in section \ref{problem_formulation}, we set up notation and terminology and formulate the model problem.
In section \ref{Steklov-metrics}, we discuss generalized Poincar\'e-Steklov operators as the basis for Riemannian metrics for shape manifolds. Section \ref{manifold} is devoted to the set of all shapes in the context of the novel metric introduced in section \ref{Steklov-metrics}. Section \ref{NLP_vol} rephrases NLP-like optimization algorithms on shape manifolds within the framework of domain integral formulations of shape derivatives. Finally, section \ref{numex} discusses algorithmic and implementation details, as well as, numerical results for a parabolic transmission shape optimization problems.

%
\section{Problem Formulation}\label{problem_formulation}

We first set up notation and terminology in shape calculus. Then we recall the model problem in \cite{Schulz-Structure-2014}, which is motivated by electrical impedance tomography and given by a parabolic interface shape optimization problem.

\subsection{Notations in shape calculus}\label{notations_definitions}
Let $d\in\N$ and $\tau>0$. We denote by $\Omega\subset \R^d$ a bounded domain with Lipschitz boundary $\Gamma:=\partial\Omega$ and by $J$ a real-valued functional depending on it.
Moreover, let $\{F_t\}_{t\in[0,\tau]}$ be a family of bijective mappings $F_t\colon \Omega\to\R^d$ such that $F_0=id$.
This family transforms the domain $\Omega$ into new perturbed domains $\Omega_t:=F_t(\Omega)=\{F_t(x)\colon x\in \Omega\}$ with $\Omega_0=\Omega$ and the boundary $\Gamma$ into new perturbed boundaries $\Gamma_t:=F_t(\Gamma)=\{F_t(x)\colon x\in \Gamma\}$ with $\Gamma_0=\Gamma$. If you consider the domain $\Omega$ as a collection of material particles, which are changing their position in the time-interval $[0,\tau]$, then the family $\{F_t\}_{t\in[0,\tau]}$ describes the motion of each particle, i.e., at the time $t\in [0,\tau]$ a material particle $x\in\Omega$ has the new position $x_t:=F_t(x)\in\Omega_t$ with $x_0=x$.
The motion of each such particle $x$ could be described by the \emph{velocity method}, i.e., as the flow $F_t(x):=\xi(t,x)$ determined by the initial value problem
\begin{equation}
\begin{split}
\frac{d\xi(t,x)}{dt}&=V(\xi(t,x))\\
\xi(0,x)&=x
\end{split}
\end{equation}
or by the \emph{perturbation of identity}, which is defined by $F_t(x):=x+tV(x)$ where $V$ denotes a sufficiently smooth vector field. We will use the perturbation of identity throughout the paper. The \emph{Eulerian derivative} of $J$ at $\Omega$ in direction $V$ is defined by
\begin{equation}
\label{eulerian}
DJ(\Omega)[V]:= \lim\limits_{t\to 0^+}\frac{J(\Omega_t)-J(\Omega)}{t}.
\end{equation}
The expression $DJ(\Omega)[V]$ is called the \emph{shape derivative} of $J$ at $\Omega$ in direction $V$ and $J$ \emph{shape differentiable} at $\Omega$ if for all directions $V$ the Eulerian derivative (\ref{eulerian}) exists and the mapping $V\mapsto DJ(\Omega)[V]$ is linear and continuous. For a thorough introduction into shape calculus, we refer to the monographs \cite{Delfour-Zolesio-2001,SokoZol}.
\begin{rev}
In particular, \cite{Sokolowski-1996} states that shape derivatives can always be expressed as boundary integrals due to the Hadamard structure theorem.
\end{rev}
The shape derivative arises in two equivalent notational forms:
\begin{align}
\label{dom_form}
DJ_\Omega[V]&:=\int_\Omega F(x)V(x)\, dx  &\text{(domain formulation)}\\
\label{bound_form}
DJ_\Gamma[V]&:=\int_\Gamma f(s)V(s)^\top n(s)\, ds  &\text{(boundary formulation)}
\end{align}
where $F(x)$ is a (differential) operator acting linearly on the perturbation vector field $V$ and $f\colon\Gamma\to\R$ with
\begin{equation}
DJ_\Omega[V]=DJ(\Omega)[V]=DJ_\Gamma[V].
\end{equation}
The boundary formulation (\ref{bound_form}), $DJ_\Gamma[V]$, acting on the normal component of $V$ has led to the interpretation as tangential vector of a corresponding shape manifold in \cite{VHS-shape-Riemann}.

\subsection{PDE model definition}
\begin{figure}[h]
\vspace*{.6cm}
\begin{center}
   \begin{overpic}[width=.5\textwidth]{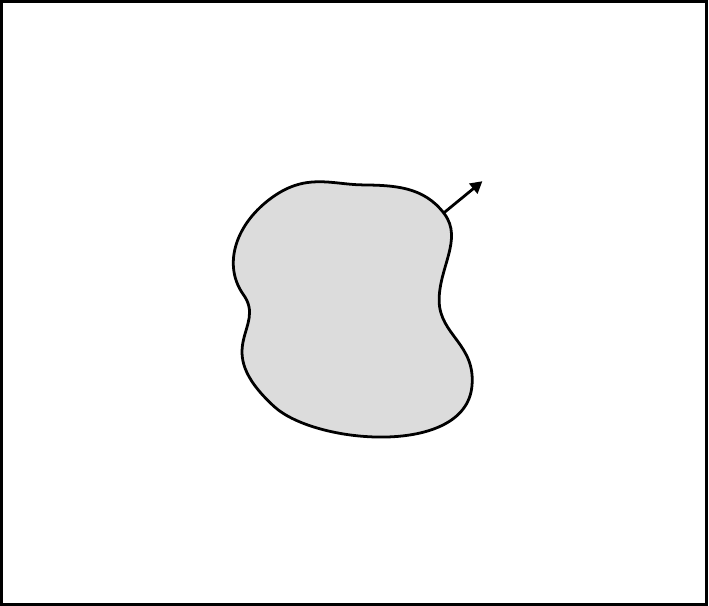}
   \put(42,40){$\Omega_2$}
   \put(20,13){$\Omega_1$}
   \put(44,88){$\Gt $}
   \put(43,61.5){$\Gi $}
   \put(-11,42){$\Gl $}
   \put(101,42){$\Gr $}
   \put(44,-6){$\Gb $}
   \put(65.5,54){$n$}
   \end{overpic}
   \end{center}
   \vspace*{.4cm}
   \caption{Example of a domain $\Omega=\Omega_1\cupdot\Gi\cupdot\Omega_2$ where $\Go:=\partial\Omega=\Gb\cupdot\Gl\cupdot\Gr\cupdot\Gt$ and $n$ denotes the unit outer normal to $\Omega_2$ at $\Gi$}
   \label{fig_Omega}
\end{figure}

We use the same model problem as in \cite{Schulz-Structure-2014}, which is briefly recalled. Let this domain $\Omega$ be an open subset of $\R^2$ and split into the two disjoint subdomains $\Omega_1,\Omega_2\subset\Omega$ such that $\Omega_1\cupdot\Gi\cupdot\Omega_2=\Omega$, $\Gb\cupdot\Gl\cupdot\Gr\cupdot\Gt=\partial\Omega\, \, (=:\Go)$ and $\partial\Omega_1\cap\partial\Omega_2=\Gi$ where the interior boundary $\Gi$ is assumed to be smooth and variable and the outer boundary $\Go$ Lipschitz and fixed. An example of such a domain is illustrated in figure \ref{fig_Omega}.
\begin{rev}
\begin{remark}
In the shape optimization method proposed in this work the topology of the domain $\Omega$ is fixed. This means we do not consider topology optimization.
\end{remark}
\end{rev}

The parabolic PDE constrained shape optimization problem is given in strong form by
\begin{align}\label{oc1}
\min \,  J(\Omega)= j(\Omega)+&j_\text{reg}(\Omega) :=\intdxdt{(y-\bar{y})^2}+\mu\int_{\Gi}1\, ds\\
\label{oc2}
\mbox{s.t. } \td{y}- \mathrm{div}(k\nabla y)&=f^{\text{model}}\quad \text{in }\Omega\times(0,T]
\\
\label{oc3}
\hspace{20mm}y&=1\quad \text{on }\Gt\times(0,T]
\\
\label{oc4}
\nd{y}&=0\quad \text{on }(\Gb\cup\Gl\cup\Gr)\times(0,T]
\\
\label{oc5}
y&=y_0\quad\text{in }\Omega\times\{0\}
\end{align}
where
\begin{equation*}
k\equiv\begin{cases}
k_1 = \mathrm{const.}\quad\text{ in }\Omega_1\times(0,T]\\
k_2 = \mathrm{const.} \quad\text{ in }\Omega_2\times(0,T]
\end{cases}
\end{equation*}
and $n$ denotes the unit outer normal to $\Omega_2$ at $\Gi$. Of course, the formulation (\ref{oc2}) of the differential equation is to be understood only formally because of the jumping coefficient $k$. We observe that the unit outer normal to $\Omega_1$ is equal to $-n$, which enables us to use only one normal $n$ for the subsequent discussions. Furthermore, we have interface conditions at the interface $\Gi$. We formulate explicitly the continuity of the state and of the flux at the boundary as
\begin{equation}
\label{oc6}
\left\llbracket y\right\rrbracket =0\, ,\quad \left\llbracket k\nd{y}\right\rrbracket =0\quad\text{on }\Gi\times(0,T]
\end{equation}
where the jump symbol $\llbracket\cdot\rrbracket$ denotes the discontinuity across the interface $\Gi$ and is defined by $\llbracket v \rrbracket := v_1-v_2$ 
where $v_1:=v\,\rule[-2mm]{.1mm}{4mm}_{\, \Omega_1}$ and $v_2:=v\,\rule[-2mm]{.1mm}{4mm}_{\, \Omega_2}$. The perimeter regularization, $j_\text{reg}(\Omega):=\mu\int_{\Gi}1\, ds$, with $\mu>0$ in the objective (\ref{oc1}) is frequently used in this kind of problems. In \cite{Sturm2013} a weaker but more complicated regularization is instrumental in order to show existence of solutions.

\begin{rev}
We assume $f^{\text{model}}\in L^2(0,T;L^2(\Omega))$ and $y\in L^2\left(0,T;H^1(\Omega)\right)$.
\end{rev}
In our setting, the boundary value problem (\ref{oc2}-\ref{oc6}) is written in weak form 
\begin{equation}
\label{wf}
a^{\text{model}}(y,p)=b^{\text{model}}(p,p^1 p^2)\, , \ \forall p\in
W\left(0,T;H^1(\Omega)\right)
\end{equation}
and for all $p^1\in W\left(0,T;H^{1/2}(\Gb\cup\Gl\cup\Gr)\right)$,
$p^2\in W\left(0,T;H^{-1/2}(\Gt)\right)$
as in \cite{Schulz-Structure-2014}. For properties of the function spaces, we refer the reader to the literature, e.g. \cite{GrossReusken,Troeltzsch}.
The bilinear form $a^{\text{model}}(y,p)$ in (\ref{wf}) is achieved by applying integration by parts on $\intdxdt{\frac{\partial y}{\partial t}p}$ and on $\intdxdt{\text{div}(k\nabla y)p}=\intdt{\int_{\Omega_1}\text{div}(k_1\nabla y_1)p_1 \ dx}+\intdt{\int_{\Omega_2}\text{div}(k_2\nabla y_2)p_2 \ dx}$.
Thus, we get
\begin{equation}
\begin{split}
a^{\text{model}}(y,p) &:=\intdx{y(T,x)\, p(T,x)}-\intdx{y_0\, p(0,x)}-\intdxdt{\td{p}y}\\
& \hspace*{.5cm}+ \intdxdt{k\nabla y^\top\nabla p}-\intdsidt{\left\llbracket k\nd{y}p\right\rrbracket}\\
&\hspace*{.5cm}-\intdsodt{k_1\nd{y}p}.\label{wfa2}
\end{split}
\end{equation}
The linear form $b^{\text{model}}(p,p^1,p^2)$ in (\ref{wf}) is given by
\begin{equation}
b^{\text{model}}(p,p^1,p^2) :=b^{\text{model}}_1(p)+b^{\text{model}}_2(p^1,p^2)
\label{wfb}
\end{equation}
where
\begin{align}
b^{\text{model}}_1(p)&:=\intdxdt{f^{\text{model}}p},\label{b1} \\
b^{\text{model}}_2(p^1,p^2)&:=\intdt{\int_{\Gt}p^1(y-1)\, ds}+\intdt{\int_{\Go\setminus\Gt}p^2\nd{y}\, ds}.\label{b2}
\end{align}
In the following, we assume for the observation \begin{rev}$\bar{y}\in L^2\left(0,T;H^1(\Omega)\right)$\end{rev}.
The Lagrangian of (\ref{oc1}-\ref{oc6}) is defined as
\begin{equation}
\LL(\Omega,y,p):=J(\Omega)+a^{\text{model}}(y,p)-b^{\text{model}}(p,p^1,p^2)
\label{lagrangian}
\end{equation}
where $J(\Omega)$ is defined in (\ref{oc1}), $a^{\text{model}}(y,p)$ in (\ref{wfa2}) and $b^{\text{model}}(p,p^1,p^2)$ in (\ref{wfb}-\ref{b2}).

The adjoint problem, which we obtain from differentiating the Lagrangian $\LL$ with respect to $y$, is given in strong form by 
\begin{align}
-\frac{\partial p}{\partial t}-\mathrm{div}(k\nabla p)&=-(y-\overline{y}) \quad \text{in }\Omega \times [0,T)\label{adjoint1}\\
p&= 0\quad\text{in }\Omega \times \{T\}\label{adjoint2}\\
\left\llbracket k\nd{p}\right\rrbracket&=0 \quad\text{on }\Gi \times [0,T)\label{adjoint3}\\
\left\llbracket p\right\rrbracket&=0 \quad\text{on }\Gi \times [0,T)\label{adjoint4}\\
\nd{p}&=0\quad\text{on }\left(\Gb\cup\Gl\cup\Gr\right) \times [0,T)\label{adjoint5}\\
p&=0\quad\text{on }\Gt \times [0,T)\label{adjoint6}\\
p^1&=-k_1p\quad\text{on }\left(\Gb\cup\Gl\cup\Gr\right) \times [0,T)\label{adjoint7}\\
p^2&=k_1\nd{p}\quad\text{on }\Gt \times [0,T)\label{adjoint8}
\end{align}
and the state equation, which we get by differentiating the Lagrangian $\LL$ with respect to $p$, is given in strong form by
\begin{equation}
\td{y}-\mathrm{div}(k\nabla y)=f^{\text{model}} \quad \text{in }\Omega \times (0,T]. \label{design}
\end{equation}

As mentioned earlier, in many cases, the shape derivative arises in two equivalent forms. If we consider the objective (\ref{oc1}) without the perimeter regularization $j_\text{reg}$, the shape derivative can be expressed as an integral over the domain $\Omega$ as well as an integral over the interface $\Gi$.
Assume that a solution $y$ of the parabolic PDE problem (\ref{oc2}-\ref{oc6}) exists and is at least in $L^2\left(0,T;H^1(\Omega)\right)$. Moreover, assume that the adjoint equation (\ref{adjoint1}-\ref{adjoint6}) admits a solution $p\in W\left(0,T;H^1(\Omega)\right)$. Then the shape derivative of the objective $J$ without perimeter regularization, i.e., the shape derivative of $j$ at $\Omega$ in the direction $V$ expressed as an integral over the domain $\Omega$ is given by
\begin{equation}
\label{boundary_expression}
\begin{split}
Dj_\Omega[V]=\int_{0}^{T}&\int_{\Omega}-k\nabla y^\top\left(\nabla V+\nabla V^\top\right)\nabla p-p\left(\nabla f^{\text{model}}\right)^\top V\\
&+\mathrm{div}(V)\left(\frac{1}{2}(y-\overline{y})^2+\td{y}p+k\nabla y^\top\nabla p-f^{\text{model}}p\right)dx\, dt.
\end{split}
\end{equation}
This domain integral allows us to calculate the boundary expression of the shape derivative, which is given by
\begin{equation}
Dj_\Gi[V]=\intdsidt{\left\llbracket k\right\rrbracket \nabla y_1^\top\nabla p_2\left<V,n\right>}.\label{shape_der2}
\end{equation}
The derivations are very technical. Note that we need a higher regularity of $y$ and $p$ to provide the boundary shape derivative expression (\ref{shape_der2}). More precisely, $p$ has to be an $L^2(0,T;H^2(\Omega))$-function having weak first derivatives in $L^2(0,T;H^1(\Omega)')$ and $y$ has to be an element of $L^2(0,T;H^2(\Omega))$. Here $H^1(\Omega)'$ denotes the dual space of $H^1(\Omega)$. We achieve (\ref{boundary_expression}) by an application of the theorem of Correa and Seeger \cite[theorem 2.1]{CorreaSeger} and (\ref{shape_der2}) by an application of integration by parts. We refer the reader for its derivations to \cite{Schulz-Structure-2014}. 
By combining theorem 2.1 and 2.2 in \cite{Schulz-Structure-2014} with proposition 5.1 in \cite{Novruzi-2002} we get the following two expressions for the shape derivative of the objective $J$ (with perimeter regularization) at $\Omega$ in the direction $V$:
\begin{equation}
\label{sd1}
Dj_\Omega[V]+Dj_\text{reg}(\Omega)[V]
=DJ(\Omega)[V]=Dj_\Gi[V]+Dj_\text{reg}(\Omega)[V]
\end{equation}
with
\begin{equation}
\label{reg1}
Dj_\text{reg}(\Omega)[V]=\int_{\Gi}\left<V,n\right>\mu\kappa\, ds
\end{equation}
where $\kappa$ denotes the curvature corresponding to the normal $n$.


%
\section{Steklov-Poincar\'e type metrics on shape manifolds}
\label{Steklov-metrics}
We first discuss the definition of shape manifolds and metrics. Then, we introduce novel metrics dovetailed to shape optimization based on domain formulations of shape derivatives.

\subsection{Shape manifolds}\label{manifold_definitions}
As pointed out in \cite{VHS-shape-Riemann}, shape optimization can be viewed as optimization on Riemannian shape manifolds and resulting optimization methods can be constructed and analyzed within this framework, which combines algorithmic ideas from \cite{Absil-book-2008} with the differential geometric point of view established in \cite{MM-2006}.
As in \cite{VHS-shape-Riemann}, we first study connected and compact subsets $\Omega_2\subset\Omega\subset\R^2$ with $\Omega_2\neq\emptyset$ and $C^\infty$ boundary $\partial\Omega_2$ where $\Omega$ denotes a bounded domain with Lipschitz-boundary (cf. figure \ref{fig_Omega}).
We now identify the variable boundary $\partial \Omega_2 = \Gi$ with a simple closed curve $c\colon S^1 \to \R^2$.
Additionally, we need to describe a space including all feasible shapes $\Gi$ and the corresponding tangent spaces.
In \cite{MM-2006}, this set of smooth boundary curves $c$ is characterized by  
\begin{equation}
B_e(S^1,\R^2):=\mbox{Emb}(S^1,\R^2)/\mbox{Diff}(S^1),
\end{equation}
i.e., as the set of all equivalence classes of $C^\infty$ embeddings of $S^1$ into the plane ($\mbox{Emb}(S^1,\R^2)$), where the equivalence relation is defined by the set of all $C^\infty$ re-parameterizations, i.e., diffeomorphisms of $S^1$ into itself ($\mbox{Diff}(S^1)$). A particular point on the manifold $B_e(S^1,\R^2)$ is represented by a curve $c\colon S^1\ni\theta\mapsto c(\theta)\in\R^2$. Because of the equivalence relation ($\mbox{Diff}(S^1)$), the tangent space is isomorphic to the set of all normal $C^\infty$ vector fields along $c$, i.e.,
\begin{equation}
\label{isomorphism_Be}
T_cB_e\cong\{h\colon h=\alpha n,\, \alpha\in C^\infty(S^1,\R)\}
\end{equation}
where $n$ is the unit exterior normal field of the shape $\Omega_2$ defined by the boundary $\partial\Omega_2=c$ such that $n (\theta)\perp c^\prime$ for all $\theta\in S^1$ and $c^\prime$ denotes the circumferential derivative as in \cite{MM-2006}. 
Several intrinsic metrics are discussed in \cite{MM-2006}, among which the following Sobolev metric seems the most natural intrinsic one from a numerical point of view. For $A>0$, the Sobolev metric is induced by the scalar product
\begin{equation}\label{eq_sobolev_metric}
\begin{split}
g^1\colon T_cB_e\times T_cB_e & \to \R,\\
    (h,k) & \mapsto ((id-A\triangle_c)\alpha,\beta)_{L^2(c)}
\end{split}
\end{equation}
where $h=\alpha n$ and $k=\beta n$ denote two elements from
the tangent space at $c$ and $\triangle_c$ denotes the Laplace-Beltrami operator on the surface $c$.  In \cite{MM-2006} it is shown that the condition $A>0$ guarantees that the scalar product $g^1$ defines a Riemannian metric on $B_e$ and 
thus, geodesics can be used to measure distances.

With the shape space $B_e$ and its tangent space in hand we can now form the Riemannian shape gradient corresponding to a shape derivative given in the form
\begin{equation}
dJ(\Omega)[V]=\int_c \psi\scp{V,n}ds.
\end{equation}
In our model setting the objective function $J$ is given in (\ref{oc1}) and its shape derivative in (\ref{sd1}).
Finally, the Riemannian shape gradient $\text{grad}J$ with respect to the Riemannian metric $g^1$ is obtained by
\begin{equation}
\grad J=gn \mbox{ with } (id -A\triangle_c)g=\psi\, .
\end{equation}
The metric $g^1$, which is also used in \cite{Schulz-Structure-2014}, necessitates a shape derivative in Hadamard form as an  efficient means to solve linear systems involving the Laplace Beltrami operator in surfaces. All of that is certainly not impossible but requires computational overhead which we can get rid of by usage of the metric discussed below. We compare the algorithmic aspects of both approaches below in section \ref{numex}.

\subsection{Steklov-Poincar\'e type Riemannian metrics}\label{manifold_metric}

The ideal Riemannian metric for shape manifolds in the context of PDE constrained shape optimization problems is to be derived from a symmetric representation of the second shape derivative in the solution of the optimization problems. Often, this operator can be related to the Dirichlet to Neumann map, aka Steklov-Poincar\'e operator, or the Laplace-Beltrami operator \cite{SS-2009}. If one aims at mesh independent convergence properties, one of these two will be appropriate in most cases. Since it can be observed that the Laplace-Beltrami operator is spectrally equivalent to the square of the Steklov-Poincar\'e operator, the latter operator seems to be more fundamental and we will focus on it as a basis for the scalar product on $T_cB_e$. Another advantage of this operator is that is blends well in with a corresponding mesh deformation strategy.

Most often, the Dirichlet to Neumann map is associated with the Laplace operator. However, as pointed out in \cite{Agoshkov-1985,KhW-2004} more general elliptic operators can be involved. For the purpose of mesh deformation, an elasticity operator may be the ideal choice. 
In numerical computations, its inverse, the Neumann to Dirichlet map or Poincar\'e-Steklov is also of importance. Therefore, we first define these operators.

In the sequel, we use the continuous generalized trace map 
\begin{equation}
\label{gamma}
\begin{split}
\gamma\colon  H^1_0(\Omega,\R^d) & \to H^{1/2}(\Gi,\R^d)\times H^{-1/2}(\Gi,\R^d),\\
U & \mapsto\mat(1){\gamma_0U\\\gamma_1U}:=\mat(1){U\,\rule[-2mm]{.1mm}{4mm}_{\, \Gi}\\\partial_nU\,\rule[-2mm]{.1mm}{4mm}_{\, \Gi}}.
\end{split}
\end{equation}

Analogously to \cite{KhW-2004}, we define for vector fields $U,V\in H^1_0(\Omega,\R^d)$ with $d=2$ or $d=3$, the Neumann solution operator for the inner boundary $\Gi$ derived from a symmetric and coercive bilinear form
\begin{equation}\label{eq_bilinearform}
a\colon H^1_0(\Omega,\R^d)\times H^1_0(\Omega,\R^d)\to\R
\end{equation}
 by
\begin{equation}
\begin{split}
E_N\colon H^{-1/2}(\Gi,\R^d) & \to H^1_0(\Omega,\R^d),\\
 u & \mapsto U
\end{split}
\end{equation}
with $U$ defined as the solution of the variational problem
\begin{equation}\label{weak-elasticity-N}
a(U,V)=\int_\Gi u^\top (\gamma_0V)\ ds\, , \ \forall\  V\in H^1_0(\Omega,\R^d)
\end{equation}
where we note that the integral in the right hand side of equation (\ref{weak-elasticity-N}) is to be understood as the duality pairing.
Furthermore, we define the Dirichlet solution operator for the inner boundary $\Gi$ by
\begin{equation}
\begin{split}
E_D\colon H^{1/2}(\Gi,\R^d) & \to H^1_0(\Omega,\R^d),\\
 u& \mapsto U
\end{split}
\end{equation}
with $U$ defined as the solution of the variational problem
\begin{equation}\label{weak-elasticity-D}
a(U,V)=0\, , \ \forall\  V\in H^1_0(\Omega,\R^d)\, ; \quad
U\,\rule[-2mm]{.1mm}{4mm}_{\, \Gi} = u.
\end{equation}

Now, we can define the Dirichlet to Neumann map and the Neumann to Dirichlet map as done in the following definition: 

\begin{definition}
In the setting above, the Dirichlet to Neumann map $T$ and the Neumann to Dirichlet map $S$ are defined by
\begin{align}
T&:=\gamma_1\circ E_D\colon H^{1/2}(\Gi,\R^d)\to H^{-1/2}(\Gi,\R^d),\\
S&:=\gamma_0\circ E_N\colon H^{-1/2}(\Gi,\R^d)\to H^{1/2}(\Gi,\R^d)
\end{align}
where $\gamma_0,\gamma_1$ are given in (\ref{gamma}).
\end{definition}

In obvious generalization of theorem 2.3.1 in \cite{KhW-2004} from scalar fields to vector fields, we conclude that both operators are symmetric w.r.t.~the standard dual pairing $\langle \cdot,\cdot\rangle$, coercive, continuous and that $T=S^{-1}$, an observation, for which \cite{Agoshkov-1985} is cited in \cite{KhW-2004}. For the purpose of defining an appropriate scalar product on the tangent space of shape manifolds, we define the following mappings.
\begin{definition}In the setting above, we define
\begin{alignat*}{3}
\eta\colon H(\Gi)& \to H(\Gi,\R^d) \qquad& \eta^\top\colon H(\Gi,\R^d) & \to H(\Gi)\\
\alpha & \mapsto \alpha\cdot n &  U & \mapsto n^\top U
\end{alignat*}
where $H\in\{H^{-1/2},H^{1/2}\}$, and thus the projected operators
\begin{align}
T^p&:=\eta^\top\circ T\circ\eta \colon H^{1/2}(\Gi)\to H^{-1/2}(\Gi),\\
S^p&:=\eta^\top\circ S\circ\eta\colon H^{-1/2}(\Gi)\to H^{1/2}(\Gi).
\end{align}
\end{definition}

Both operators, $T^p$ and $S^p$, inherit symmetry, coercivity, continuity and invertibility from the operators $T,S$. However, we observe in general $T^p\neq(S^p)^{-1}$. Both operators can be used for the definition of a scalar product on the tangent space. In line with the discussion of Sobolev type metrics in \cite{MM-2006}, we would prefer a scalar product with a smoothing effect like the projected Dirichlet to Neumann map $T^p$. However, we need its inverse in numerical computations, which is usually not $S^p$, although spectrally equivalent. We can limit the computational burden, if we use directly $(S^p)^{-1}$ as a metric on the tangent space, having a similar smoothing effect but also the advantage of the straight forward inverse $S^p$. In order to summarize, let us explicitly formulate the operator
\begin{equation}
\begin{split}
S^p\colon  H^{-1/2}(\Gi)& \to H^{1/2}(\Gi),\\
         \alpha & \mapsto (\gamma_0 U)^\top n
\end{split}
\end{equation}
where $U\in H^1_0(\Omega,\R^d)$ solves the Neumann problem
\begin{equation}\label{weak-elasticity-N2}
a(U,V)=\int_\Gi \alpha\cdot (\gamma_0 V)^\top n\ ds\, , \ \forall\  V\in H^1_0(\Omega,\R^d)
\end{equation}
which corresponds to an elliptic problem with fixed outer boundary and forces $\alpha\cdot n $ at the inner boundary $\Gi$. 
Thus, we propose to use the scalar product $g^S$ defined below.
\begin{definition} In the setting above, we define the scalar product $g^S$ on $H^{1/2}(\Gi)$ by
\begin{equation}\label{scp}
\begin{split}
g^S\colon H^{1/2}(\Gi)\times H^{1/2}(\Gi) & \to \R,\\
(\alpha,\beta) &\mapsto \langle\alpha,(S^p)^{-1}\beta\rangle=
\int_\Gi \alpha(s)\cdot [(S^p)^{-1}\beta](s)\ ds.
\end{split}
\end{equation}
\end{definition}



%
\section{Shape quasi-Newton methods based on the metric $g^S$}
\label{NLP_vol}

\begin{rev}
As already mentioned in section \ref{problem_formulation} the shape derivative can always be expressed as boundary integral
$DJ_\Gi[V]=\int_{\Gi}f\left<V, n\right>  ds$ (cf. (\ref{bound_form}))
due to the Hadamard structure theorem. If $V\,\rule[-2mm]{.1mm}{4mm}_{\hspace{.6mm}\Gi}=\alpha n$, this can be written more concisely as
\begin{equation}
\label{HadamardConcisely}
DJ_{\Gi}[V]=\int_{\Gi} \alpha f \ ds.
\end{equation}
Due to the isomorphism (\ref{isomorphism_Be}) and the handy expression (\ref{HadamardConcisely}) we can state the connection of $\left(B_e(S^1,\R^2),g^S\right)$ with shape calculus, i.e., we can determine a representation $h\in T_{\Gi}B_e(S^1,\R^2) \cong \{h\colon h=\alpha n,\, \alpha\in C^\infty(S^1,\R)\}$ of the shape gradient in terms of $g^S$ defined in (\ref{scp}) by
\begin{equation}
g^S(\phi,h)=(f,\phi)_{L^2(\Gi)},\quad \forall \phi\in C^\infty(\Gi,\mathbb{R}),
\end{equation}
which is equivalent to
\begin{equation}
\int_\Gi \phi(s)\cdot [(S^p)^{-1}h](s) \ ds=\int_\Gi f(s)\phi(s) \ ds ,\quad \forall \phi\in C^\infty(\Gi,\mathbb{R}).
\end{equation}
\end{rev}
Thus, $h=S^pf=(\gamma_0 U)^\top n$, where $U\in H^1_0(\Omega,\R^d)$ solves
\begin{equation}\label{eq_most_important}
a(U,V)=\int_\Gi f\cdot (\gamma_0 V)^\top n \ ds
=DJ_\Gi[V]=DJ_\Omega[V]\, , \ \forall\  V\in H^1_0(\Omega,\R^d)
\end{equation}
which means that the representation of the domain integral formulation in terms of the elliptic form $a(\cdot,\cdot)$ as used in \cite{Langer-2015} can be -- projected to the normal component on $\Gi$ -- interpreted as the representation of the boundary integral formulation in terms of $(S^p)^{-1}$. However, in both points of view, the information of the shape derivative is in physical terms used as a force (in the domain or on the boundary) and we obtain a vector field $U$ as an (intermediate) result, which can serve as a deformation of the computational mesh -- identical to Dirichlet deformation.

\begin{rev}
\begin{remark}\label{remark_injectivity}
In general, $h=S^pf=(\gamma_0 U)^\top n$ is not necessarily an element of $T_{\Gi}B_e$ because it is not ensured that $U\in H^1_0(\Omega,\mathbb{R}^d)$ is $C^\infty$. Under special assumptions depending on the coefficients of a second-order partial differential operator, the right hand-side of a PDE and the domain $\Omega$ on which the PDE is defined, a weak solution $U\in H^1_0(\Omega,\mathbb{R}^d)$ of a PDE is $C^\infty$ by the theorem of infinite differentiability up to the boundary \cite[theorem 6, section 6.3]{Evan}.
\end{remark}
\end{rev}

Now, we rephrase the l-BFGS-quasi-Newton method for shape optimization from \cite{Schulz-Structure-2014} in terms of the metric $g^S$ and in generalization to domain formulations of the shape derivative. We note that the complete deformation of a shape optimization algorithms is just the (linear) sum of all iterations, which means that the BFGS update formulas can be rephrased directly in terms of the deformation vector field, rather than only as boundary deformations to be transferred to the domain mesh in each iteration. 

BFGS update formulas need the evaluation of scalar products, where at least one argument is a gradient-type vector. According to the metric introduced in section \ref{Steklov-metrics}, we can assume that a gradient type vector $u\in T_{c}B_e$ can be written as
\begin{equation}\label{grad-type}
u = (\gamma_0U)^\top n
\end{equation}
for some vector field $U\in H^1_0(\Omega,\R^d)$. The other argument $v$ is either of gradient-type or deformation-type, which can also be assumed of being of the form (\ref{grad-type}), i.e., 
\begin{equation}
v = (\gamma_0V)^\top n
\end{equation} 
for some $V\in H^1_0(\Omega,\R^d)$. If $u$ is a gradient of a shape objective $J$, we observe
\begin{equation}
g^S(u,v)=DJ_\Gi[V]=DJ_\Omega[V]=a(U,V).
\end{equation}
This observation can be used to reformulate the scalar product $g^S(\cdot,\cdot)$ on the boundary equivalently as $a(\cdot,\cdot)$ for domain representations. In the sequel, we consider only domain representations $U_j\in H^1_0(\Omega,\R^d)$ of $\grad J(c_j)\in H^{1/2}(\Gi)$, mesh deformations $S_j\in H^1_0(\Omega,\R^d)$ and differences $Y_j:=U_{j+1}-{\cal T}_{S_{j}}U_j\in H^1_0(\Omega,\R^d)$ where ${\cal T}_{S_{j}}$ denotes the vector transport as in \cite{Schulz-Structure-2014}.

With this notation we formulate the double-loop of an l-BFGS quasi-Newton method:

\begin{algorithmic}
\State $\rho_j \gets g^S\left((\gamma_0Y_j)^\top n,(\gamma_0S_j)^\top n\right)^{-1}=a(Y_j,S_j)^{-1}$
\State $q \gets U_j$
\For{$i = j-1, \dots , j-m$}
	\State $S_i \gets {\cal T}_{q}S_i$
	\State $Y_i \gets {\cal T}_{q}Y_i$
	\State $\alpha_i \gets \rho_i  g^S\left((\gamma_0S_i)^\top n,(\gamma_0q)^\top n\right)=\rho_i  a(S_i,q)$
	\State $q \gets q - \alpha_i Y_i$
\EndFor
\State $q \gets \frac{g^S\left((\gamma_0Y_{j-1})^\top n,(\gamma_0S_{j-1})^\top n\right)}{g^S\left((\gamma_0Y_{j-1})^\top n,(\gamma_0Y_{j-1})^\top n\right)} \ q
=\frac{a(Y_{j-1},S_{j-1})}{a(Y_{j-1},Y_{j-1})} \ q$
\For{$i = j-m, \dots , j-1$}
	\State $\beta_i \gets \rho_i g^S\left((\gamma_0Y_{i})^\top n,(\gamma_0z)^\top n\right)=\rho_i a(Y_{i},q)$
	\State $q \gets q + (\alpha_i - \beta_i) Y_i$
\EndFor \\
\Return $q = G_j^{-1} \grad J(c_j)$
\end{algorithmic}
The resulting vector $q\in H^1_0(\Omega,\R^d)$ is simultaneously a shape deformation as well as a deformation of the domain mesh.

%
\section{Numerical results and implementation details}\label{numex}

We compare the limited memory BFGS shape optimization algorithms of \cite{Schulz-Structure-2014} with the analogous algorithm based on the Riemannian metric $g^S$, introduced above. We use a test case within the domain
$\Omega=(-1,1)^2$, which contains a compact and closed subset $\Omega_2$ with smooth boundary. The parameter $k_1$ is valid in the exterior $\Omega_1 = \Omega \setminus \Omega_2$ and the parameter $k_2$ is valid in the interior $\Omega_2$. First, we build artificial data $\bar{y}$, by solving the state equation for the setting
$\bar{\Omega}_2:=\{x\colon \norm{x}_2 \le r\}$ with $r=0.5$. Afterwards, we choose another initial domain $\Omega_1$ and $\Omega_2$. Figure \ref{fig_Omega} illustrates the interior boundary $\Gamma_\text{int}$ around the initial domain $\Omega_2$ and the target domain $\bar{\Omega}_2$.
The reason for this choice of artificial test data is that we obtain a representation of $\bar{y}$ that can be evaluated at arbitrary points in space since it is represented in finite element basis.
Moreover, this construction guarantees that the optimization converges to a reasonable shape that is within the boundaries $\Omega$ and not too different to the initial shape such that the mesh remains feasible under deformations.

\begin{rev}
\begin{remark}
Choosing the measurements $\bar{y}$ as the solution of the model equation (\ref{oc2}-\ref{oc6}) we obtain that $\bar{y} \in L^2(0,T; H^1(\Omega))$ as we assumed in section \ref{notations_definitions}.
\end{remark}
\end{rev}

\begin{figure}
\begin{subfigure}{1.0\textwidth}
\begin{center}
\begin{tabular}{cccc}
\includegraphics[width=0.23\textwidth]{nz0000}&
\includegraphics[width=0.23\textwidth]{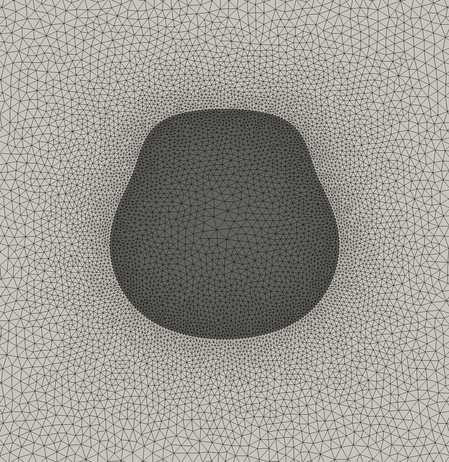}&
\includegraphics[width=0.23\textwidth]{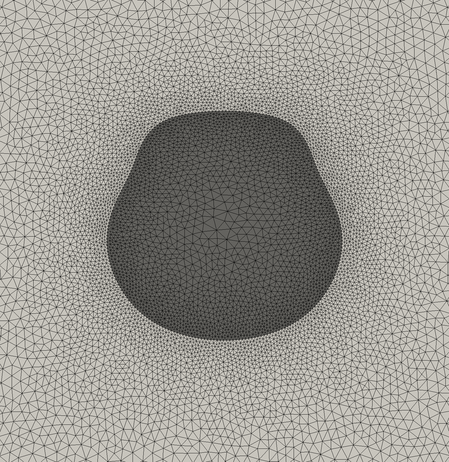}&
\includegraphics[width=0.23\textwidth]{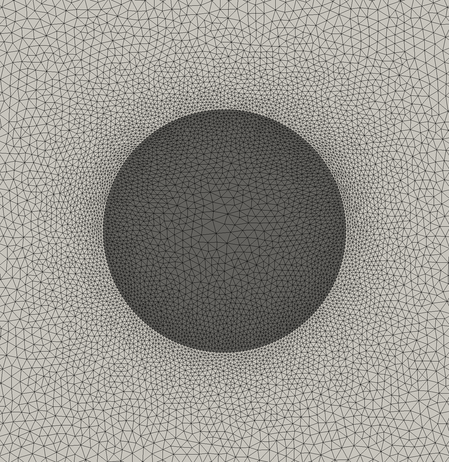}
\end{tabular}
\end{center}
\subcaption{Smooth deformations and convergence to optimal shape due to modified source term}
\label{fig_bfgs_deformations}
\end{subfigure}
\begin{subfigure}{1.0\textwidth}
\begin{center}
\begin{tabular}{cccc}
\includegraphics[width=0.23\textwidth]{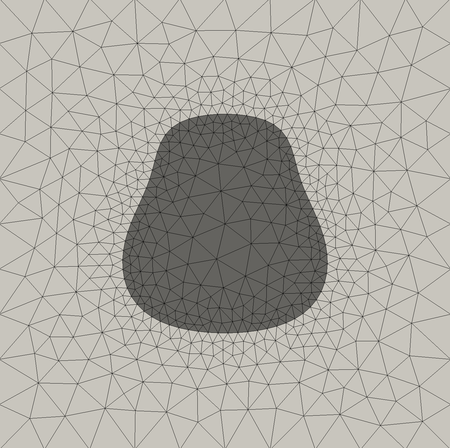}&
\includegraphics[width=0.23\textwidth]{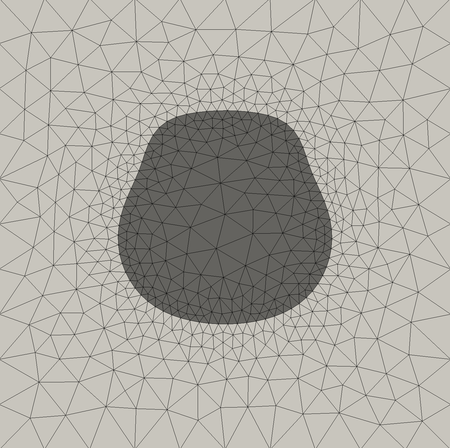}&
\includegraphics[width=0.23\textwidth]{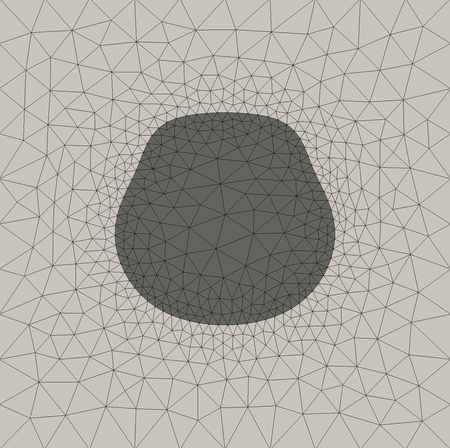}&
\includegraphics[width=0.23\textwidth]{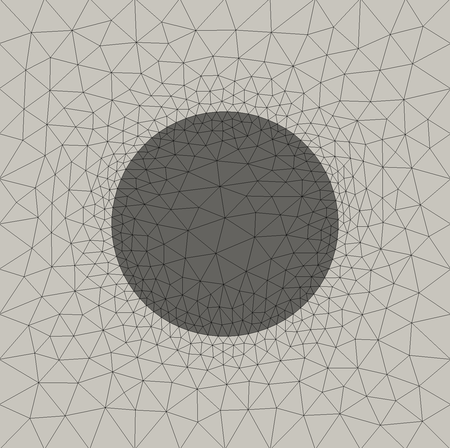}
\end{tabular}
\end{center}
\subcaption{The domain-based optimization approach also enables the use of much coarser spatial discretizations (approx.\ 1000 cells)}
\label{fig_coarse_deformations}
\end{subfigure}
\caption{BFGS iterates with corrected source term \eqref{boundary_expression} indicating mesh independent convergence}
\label{fig_meshes}
\end{figure}

\begin{figure}
\begin{center}
\begin{tabular}{cccc}
\includegraphics[width=0.23\textwidth]{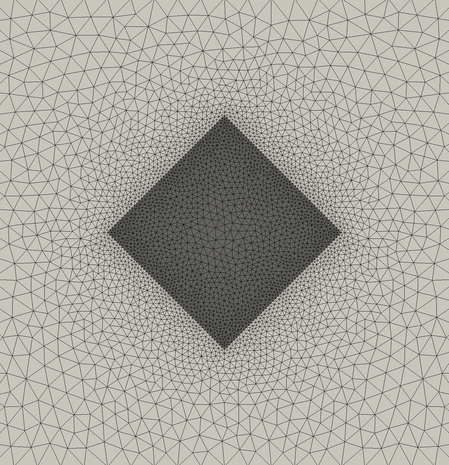}&
\includegraphics[width=0.23\textwidth]{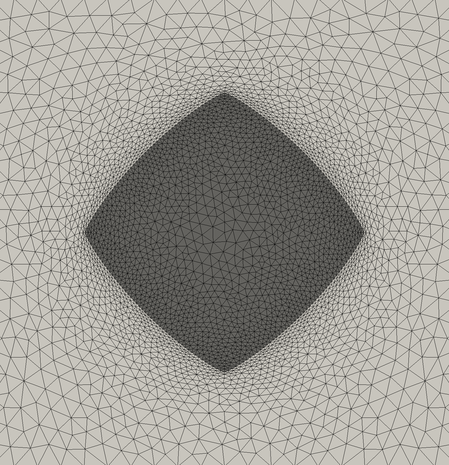}&
\includegraphics[width=0.23\textwidth]{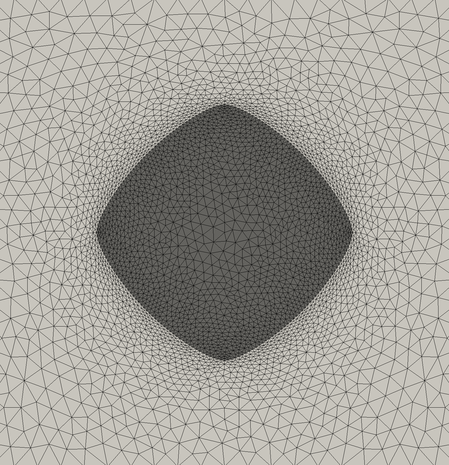}&
\includegraphics[width=0.23\textwidth]{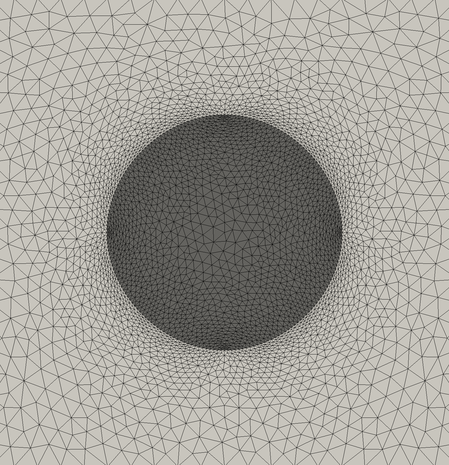}
\end{tabular}
\end{center}
\caption{Smooth mesh deformations even with kinks in the initial configuration due to regularization}
\label{fig_meshes_with_kinks}
\end{figure}

In the particular test case, which is studied in this section and can be seen in figure \ref{fig_meshes}, the diffusion coefficients are chosen to be $k_1 = 1$ and $k_2 = 0.001$.
Further, the initial condition is $y_0(x) = 0$ for all $x\in \Omega$, $f^\text{model}(x,t) = 0$ in $ (x,t) \in \Omega\times(0,T]$ and the final time of the simulation is $T = 20$.
The results shown in this section are computed under a mild perimeter regularization with $\mu = 10^{-6}$, where we did not notice any numerical difference with the case $\mu=0$.
Yet, for the non-smooth initial configuration shown in figure \ref{fig_meshes_with_kinks} a stronger regularization has to be chosen in the first iterations as $\mu_\text{init} = 0.01$.
In this particular case the regularization is controlled by a decreasing sequence from $\mu_\text{init}$ to $\mu$.

The numerical solution of the boundary value problem (\ref{oc2}-\ref{oc5}) is obtained by discretizing its weak formulation \eqref{wf} with linear finite elements in space and an implicit Euler scheme in time.
For the time discretization 30 time steps are chosen, which are equidistantly distributed.
The diffusion parameter $k$ is discretized as a piecewise constant function in contrast to the continuous trial and test functions.
This choice of function spaces ensures that the transmission conditions \eqref{oc6} are automatically fulfilled.
The corresponding adjoint problem (\ref{adjoint1}-\ref{adjoint8}) can be discretized in the same way.
More precisely, it is not necessary to assemble different linear operators, which is attractive in terms of computational effort.
All arising linear systems are then solved using the preconditioned conjugate gradient method.

\begin{rev}
Our investigations focus on the comparison between two l-BFGS optimization approaches:
The \textbf{first approach} is based on the surface expression of the shape derivative, as intensively described in \cite{Schulz-Structure-2014}.
Here, a representation of the shape gradient at $\Gi$ with respect to the Sobolev metric \eqref{eq_sobolev_metric} is computed and applied as a Dirichlet boundary condition in the linear elasticity mesh deformation.
This involves two operations, which are non-standard in finite element tools and thus leads to additional coding effort.
Since we are dealing with linear finite elements the gradient expressions of state $y$ and adjoint $p$ in \eqref{shape_der2} are piecewise constant and can not be applied directly to the mesh as deformations.
We thus have to implement a kind of $L^2$-projection on $\Gi$ (cf.~\cite{Schulz-Structure-2014}) bringing back the sensitivity information into the space of continuous, linear functions.
The next additional piece of code is a discrete version of the Laplace-Beltrami operator for the Sobolev metric \eqref{eq_sobolev_metric}.
The essential part of this is the solution of a tangential Laplace equation on the surface $\Gi$.
Therefore, we follow the presentations \cite{meyer2003discrete} and artificially extend our 2D grid in the third coordinate direction.
The \textbf{second approach}, discussed in sections \ref{Steklov-metrics} and \ref{NLP_vol}, involves the volume formulation of the shape derivative and a corresponding metric, which is very attractive from a computational point of view.
The computation of a representation of the shape gradient with respect to the chosen inner product of the tangent space is now moved into the mesh deformation itself.
The elliptic operator $a(\cdot,\cdot)$ (cf.\ \eqref{eq_bilinearform}) -- here the linear elasticity -- is both used as inner product and mesh deformation leading to only one linear system, which has to be solved.
Besides saving brain work in the calculation of the shape derivative, a lot of coding work is obsolete using surface formulation of shape derivatives.
Moreover, it is not always clear how the surface formulation looks like and which additional assumptions have to be made in its derivation.
A discussion of the l-BFGS algorithm used within this algorithm can be seen in section \ref{NLP_vol}.
\end{rev}

An essential part of a shape optimization algorithm is to update the finite element mesh after each iteration.
For this purpose, we use a solution of the linear elasticity equation
\begin{equation}\label{eq_linelas}
\begin{aligned}
\text{div}( \sigma ) &= f^\text{elas} \quad \text{in} \quad \Omega\\
U &= 0 \quad \text{on} \quad \Go\\[10pt]
\sigma &:= \lambda \text{Tr}(\epsilon) I + 2 \mu \epsilon\\
\epsilon &:= \frac{1}{2}(\nabla U + \nabla U^T)\\
\end{aligned}
\end{equation}
where $\sigma$ and $\epsilon$ are the strain and stress tensor, respectively.
Here $\lambda$ and $\mu$ denote the Lam\'{e} parameters, which can be expressed
in terms of Young's modulus $E$ and Poisson's ratio $\nu$ as
\begin{equation}
\lambda = \frac{\nu E}{(1+\nu)(1-2\nu)} \, , \quad \mu = \frac{E}{2(1+\nu)}.
\end{equation}
The solution $U\colon \Omega \to \mathbbm{R}^3$ is then added to the coordinates of the finite element nodes.
The Lam\'{e} parameters do not need to have a physical meaning here.
It is rather essential to understand their effect on the mesh deformation.
$E$ states the stiffness of the material, which enables to control the step size for the shape update.
$\nu$ gives the ratio how much the mesh expands in the remaining coordinate directions when compressed in one particular direction.
The numerical results in this work are obtained using $\nu = 0.01$ and $E = 0.1$.

Equation \eqref{eq_linelas} is modified according to the optimization approach under consideration.
In case we use the surface formulation of the shape derivative (\ref{shape_der2}), the following Dirichlet condition is added on the variable boundary
\begin{equation}\label{eq_deform_dirichlet}
U = U^\text{surf} \quad \text{on} \quad \Gi
\end{equation}
where $U^\text{surf}$ is the representation of the shape gradient with respect to the Sobolev metric $g^1$ given in \eqref{eq_sobolev_metric}.
The source term $f^\text{elas}$ is then set to zero.
Otherwise, when the mesh deformation operator is also used as shape metric, $f^\text{elas}$ assembled according to \eqref{boundary_expression} and there is no Dirichlet condition on $U$.
This only covers the portion of the shape derivative for which a volume formulation is available.
Parts of the objective function leading only to surface expressions, such as, for instance, the perimeter regularization $j_\text{reg}$, are incorporated as Neumann boundary conditions given by
\begin{equation}
\frac{\partial U}{\partial n} = f^\text{surf} \quad \text{on} \quad \Gi.
\end{equation}
\begin{rev}
In the notation of section \ref{manifold_metric} we set $a(\cdot,\cdot)$ as the weak form of the linear elasticity equation leading to
\begin{equation}
a(U,V) = \intdx{\sigma(U) : \epsilon(V)}.
\end{equation}
For our model problem given in section \ref{problem_formulation} we have to solve in the context of the domain formulation of the shape derivative and its representation in terms of $g^S$
\begin{equation}\label{eq_mesh_deform_weak}
a(U,V) = Dj_\Omega[V]+Dj_\text{reg}(\Omega)[V] \, , \ \forall\  V\in H^1_0(\Omega,\R^d)
\end{equation}
where the right hand side is given by the left formulation in \eqref{sd1}, in particular $Dj_\Omega[V]$ in (\ref{boundary_expression}) and $Dj_\text{reg}(\Omega)[V]$ in (\ref{reg1}).
Note that \eqref{eq_mesh_deform_weak} is justified by the main result \eqref{eq_most_important} of the previous sections stating that the connection between the volume formulation of shape derivatives and a bilinear form $a$ leads to a representation of the shape gradient with respect to $g^S$.
\end{rev}

\begin{figure}
\setlength{\tabcolsep}{2pt}
\begin{subfigure}{1.0\textwidth}
\begin{center}
\begin{tabular}{cccc}
\includegraphics[width=0.23\linewidth]{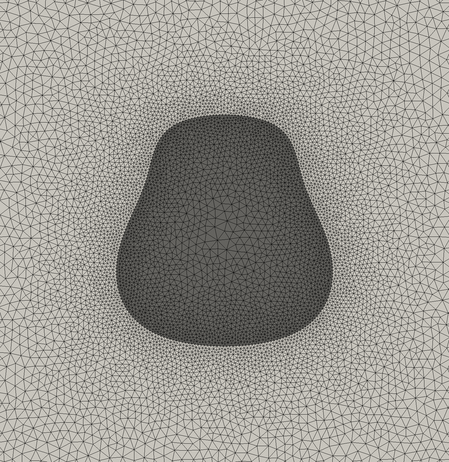}&
\includegraphics[width=0.23\linewidth]{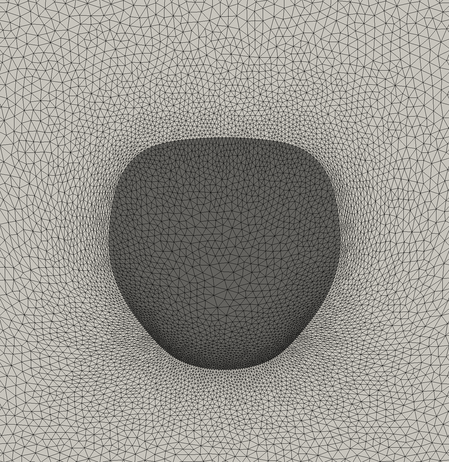}&
\includegraphics[width=0.23\linewidth]{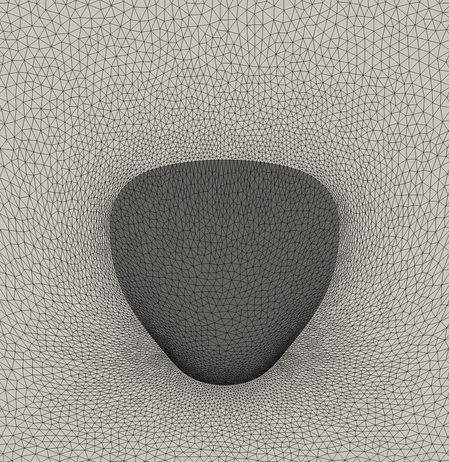}&
\includegraphics[width=0.23\linewidth]{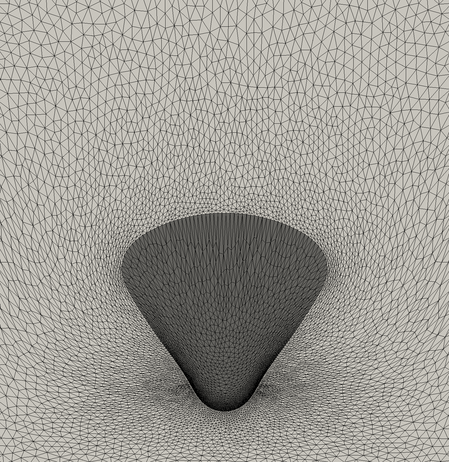}
\end{tabular}
\end{center}
\subcaption{BFGS iterates with unmodified approximation of volume shape derivative}
\label{fig_wrong_deformations}
\end{subfigure}\\
\begin{subfigure}{1.0\textwidth}
\begin{center}
\begin{tabular}{cccc}
\includegraphics[width=0.23\linewidth]{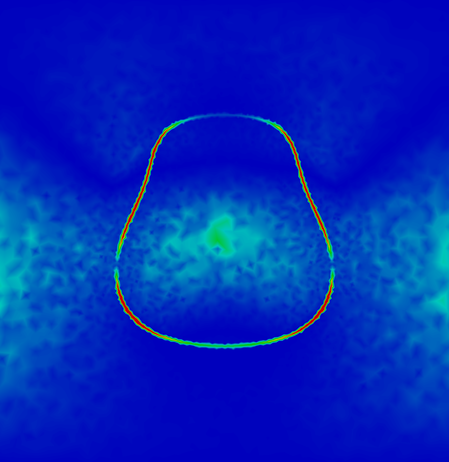}&
\includegraphics[width=0.23\linewidth]{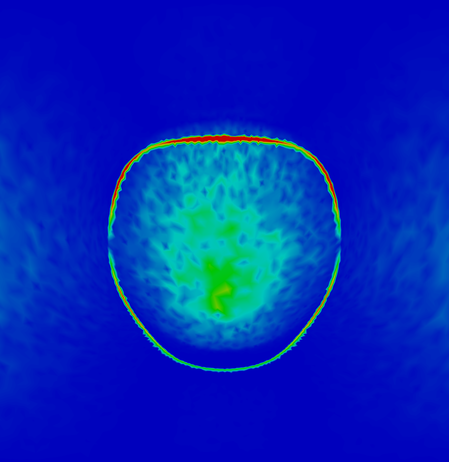}&
\includegraphics[width=0.23\linewidth]{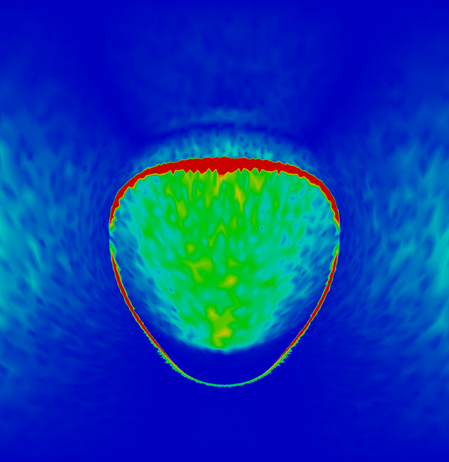}&
\includegraphics[width=0.23\linewidth]{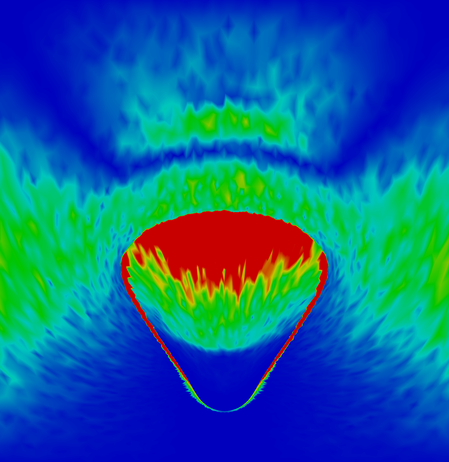}
\end{tabular}
\end{center}
\subcaption{Magnitude of unmodified volumic source term}
\label{fig_wrong_source_term}
\end{subfigure}
\caption{Wrong mesh deformations and source term due to discretization errors in the unmodified right hand side of \eqref{eq_mesh_deform_weak}}
\label{fig_wrong_meshes}
\end{figure}

Both approaches (A versus B below) follow roughly the same steps with a major difference in the way the shape sensitivity is incorporated into the mesh deformation. The appraoch A (domain formulation) is clearly to be preferred because of its implementational ease and less computational effort, if a technical detail discussed below is taken into account.
One optimization iteration can be summarized as follows:
\medskip
\begin{enumerate}
\itemsep4pt
\item[1.] The measured data $\bar{y}(t,x)$ has to be interpolated to the current iterated mesh and the corresponding finite element space.
Here, this consists of the interpolation between two finite element spaces on non-matching grids.
\item[2.] The state and the adjoint equation are solved.
\item[3.] Assembly of the linear elasticity equation.\\[3pt]
A) Domain formulation: 
\vspace{2pt}
\begin{itemize}
\itemsep2pt
\item The volume form of the shape derivative is assembled into a source term for the linear elasticity mesh deformation. Only test functions whose support includes $\Gi$ are considered, which is justified in the subsequent discussion.
The behavior of the algorithm with full assembly for all test functions is illustrated in figure \ref{fig_wrong_source_term}. Here, the magnitude of the unmodified discretization of the source term is visualized, which shows not only non-zero values outside of $\Gi$ due to discretization errors, but leads also to detrimental mesh deformations.
\item Shape derivative contributions, which are only available in surface formulation, such as the perimeter regularization, are assembled into the right hand side in form of Neumann boundary conditions.
\end{itemize}
\vspace{3pt}
B) Surface formulation:
\vspace{2pt}
\begin{itemize}
\itemsep2pt
\item The preliminary gradient $\tilde{s} = \nabla y_1^T\nabla p_2 n$ given in \eqref{shape_der2} is evaluated at $\Gi$.
\item The $L^2$-projection of $\tilde{s}$ into the space of piecewise linear, continuous functions is conducted. Let this be denoted by $\hat{s}$.
\item Finally the contributions resulting from the regularization, which is here $\kappa n$, is added to $\hat{s}$ and we solve the Laplace-Beltrami equation $(\text{id} - A\triangle_c) U^\text{surf} = \hat{s} + \kappa n$ to obtain a representation of the gradient with respect to the Sobolev metric as given in \eqref{eq_sobolev_metric}.
\item $U^\text{surf}$ then yields the Dirichlet boundary condition \eqref{eq_deform_dirichlet}.
\end{itemize}
\item[4.] Solve linear elasticity equations, apply the resulting deformation to the current finite element mesh and go to the next iteration.
\end{enumerate} 
\medskip
%

Assembling the right hand side of the discretized weak form (equation \eqref{eq_mesh_deform_weak}) only for test functions whose support intersects with $\Gi$ in the volume formulation of step 3 above is due to the following reasoning. In exact integration, the integral $Dj_\Omega[V]$ should be zero for all test functions $V$ which do not have $\Gi$ within their support. Thus, nonzero integral contributions are caused by discretization noise, On the other hand, its effect on the optimization algorithm can be understood from a perturbation point of view. We may assume that the Riemannian shape Hessian $\nabla^c \grad J$ (where $\nabla^c$ means covariant derivative), whose action in the optimal solution coincides with the action of the shape Hessian, i.e.,
\begin{equation}
g^S(\nabla^c \grad J[V],U)=D(DJ[V])[U]
\end{equation}
is coercive on the boundary, i.e., for projections $\eta^\top V|_\Gi,\eta^\top V|_\Gi$, which guarantees a well-posed problem. However, the Hessian operator approximated in the BFGS update strategy described in section \ref{NLP_vol} deals with a Hessian defined on the whole mesh, which posseses a huge kernel, determined by all vector fields with zero normal component on the boundary. Thus, the space $H^1_0(\Omega, \R^d)$ of all admissible deformations has a decomposition
\begin{equation}
H^1_0(\Omega,\R^d)=H_\Gi\oplus H_\Gi^{\perp}
\end{equation}
where
$H_\Gi:= \{E_N(\alpha n)\colon \alpha\in H^{-1/2}(\Gi)\}$
and $H_\Gi^{\perp}$ denotes its orthogonal complement in the bilinear form $a(\cdot,\cdot)$.
Shape gradients and increments in $H^1_0(\Omega,\R^d)$ lie in $H_\Gi$ only.
It is abvious that l-BFGS update formulas produce steps which lie again in $H_\Gi$ only, which means that the optimization algorithm in function spaces acts always on the coercive shape Hessian only. 
However, the discretized version is a perturbation of the infinite Hessian. Thus, perturbed coercive operators stay coercive, if the perturbation is not too large. But, positive semidefinite operators with a nontrivial kernel, almost inevitably will get directions of negative curvature, when perturbed. These directions of negative curvature will be chosen, if we allow nonzero components in the right hand side of the discretized mesh deformation equation \eqref{eq_mesh_deform_weak} in the interior of the domain. On the other hand, if we do not allow zero components there, the algorithm only acts in the subspace of the discretization of $H_\Gi$ where the projected Hessian is a perturbation of the shape Hessian and thus coercive, if the perturbation is not too large.  

\begin{figure}
\begin{center}
\includegraphics[width=1.0\linewidth]{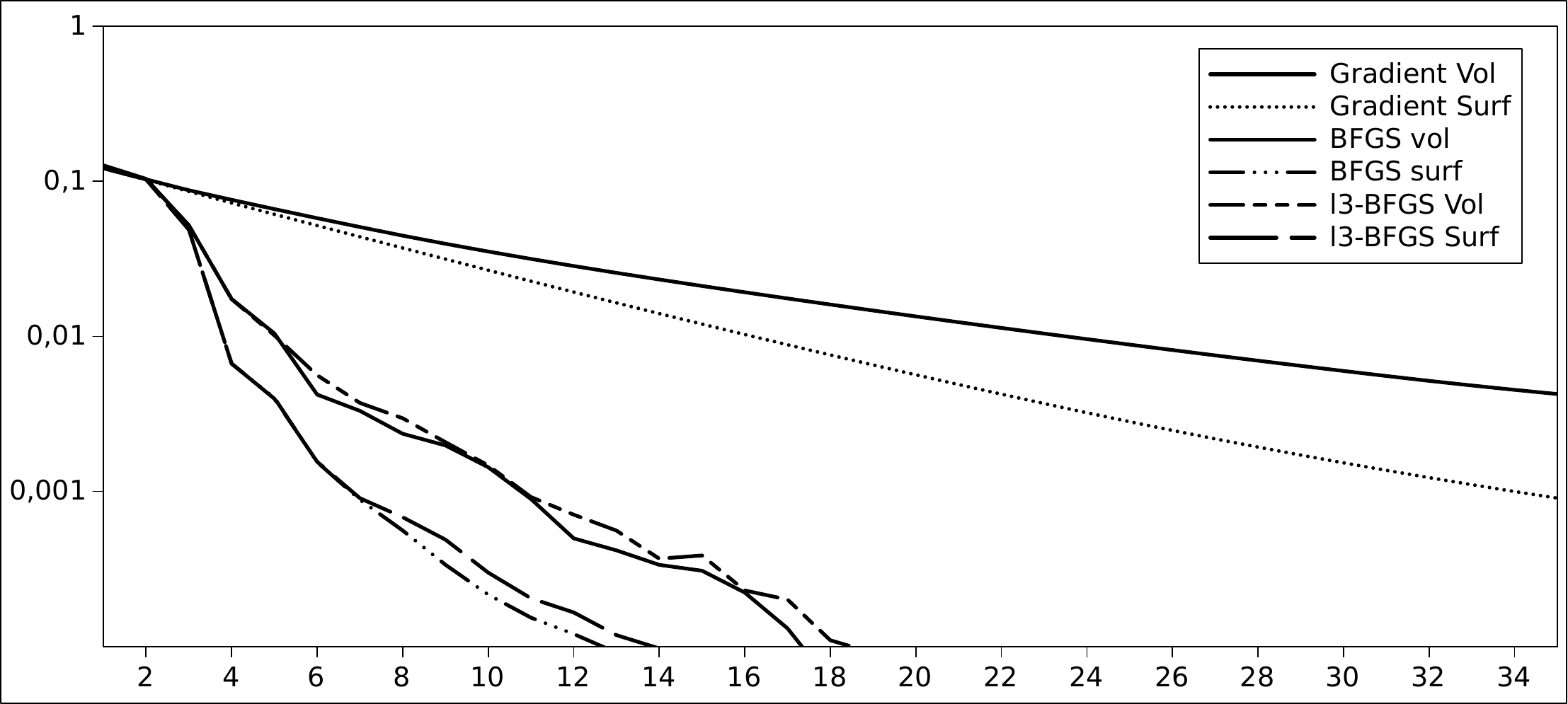}
\end{center}
\caption{Convergence of the optimization iteration measured as an approximation to the geodesic distance in the shape space on a grid with approx.\ 100,000 cells}
\label{fig_convergence_parabolic}
\end{figure}

We now conclude this section with a discussion of the numerical results.
The figures \ref{fig_wrong_meshes} to \ref{fig_meshes_with_kinks} show the initial configuration and the iterations 2, 4, 20 of the full BFGS algorithm as described in section \ref{NLP_vol}.
In figure \ref{fig_wrong_meshes} the algorithm is shown for the unmodified assembly of the right hand side in \eqref{eq_mesh_deform_weak} leading to divergence.
Whereas, figure \ref{fig_meshes} shows a selection of BFGS iterates for the modified source term and with $\mu=0$.
It is also demonstrated here that the domain-based shape optimization algorithm can be applied to very coarse meshes.
This is due to the fact, that there is no dependence on normal vectors like in the case of surface shape derivatives.
Finally, figure \ref{fig_convergence_parabolic} shows the convergence of l-BFGS with three gradients in memory, full BFGS and the pure gradient method for suface and volume shape derivative formulation, respectively.

In our tests, the convergence with the Laplace-Beltrami representation of the shape gradient seems to require a bit less iterations compared to the domain-based formulation.
Yet, the domain-based form is computationally more attractive since it also works for much coarser discretizations.
This can be seen in figure \ref{fig_meshes} where \ref{fig_bfgs_deformations} shows the necessary fineness of the mesh for the surface derivative to lead to a reasonable convergence.
The coarse grid in \ref{fig_coarse_deformations}, however, only works for the domain-based formulation.

\begin{rev}
Since the volume term $Dj_\Omega[V]+Dj_\text{reg}(\Omega)[V]$ in approach A only has to be computed for discretization elements adjacent to $\Gi$ it is computationally not more expensive than the surface formulation in approach B.
Moreover, the computing time for $L^2$-projection and solution of tangential Laplace equation is saved.
Yet, the BFGS update algorithm in approach A is more expensive then the one in B due to the higher dimension of the involved matrix, which does not play a decisive role since we only store a few gradients in memory.
These differences should yet not be overrated.
The most expensive operation is the computation of the mesh deformation involving the solution of the linear elasticity equation in both approaches A and B making them comparable in terms of computational costs.
\end{rev}

This changes for highly parallel application on supercomputers as investigated in \cite{neagel20015scalable}.
Operations, which are only performed on surfaces, can drastically affect the scalability of the overall algorithm if the computational load is not balanced also with respect to surface elements.
The higher demand for memory of the domain-based formulation seems also not be dramatic since the numerical tests suggest that very few gradients in memory are sufficient for good performance of the l-BFGS method (see figure \ref{fig_convergence_parabolic}).


%
\section{Towards a novel shape space}\label{manifold}
\begin{rev}
The scalar product introduced above in section \ref{Steklov-metrics} connects shape gradients with $H^1$ deformations. These deformations evaluated at a prior shape $\Gamma_0$ give deformed shapes $\Gi$ of class $H^{1/2}$, if the deformations are injective and continuous. In the following, it is clarified what we mean by $H^{1/2}$-shapes. The investigations done in the previous section are not limited to $C^\infty$ shapes, i.e., elements of the shape space $B_e(S^1,\mathbb{R}^2)$. Therefore, this section is devoted to an extension of $B_e(S^1,\mathbb{R}^2)$, i.e., to a novel shape space definition, and its connection to shape calculus.

We would like to recall once again that a shape in the sense of the shape space of Peter W.~Michor and David Mumford introduced in \cite{MM-2006} is given by the image of an embedding from the unit sphere $S^{d-1}$ into the Euclidean space $\R^d$. In view of our generalization, it has technical advantages to consider a prior shape $\Gamma_0$ as the boundary $\Gamma_0=\partial\mathcal{X}_0$ of a connected and compact subset $\mathcal{X}_0\subset \Omega\subset\R^d$ with $\mathcal{X}_0\neq\emptyset$, where $\Omega$ denotes a bounded Lipschitz domain. 
Let the prior set $\mathcal{X}_0$ be a Lipschitz domain, i.e., $\Gamma_0$ is a Lipschitz boundary.
An example of a prior shape is the cube. It is the union of six faces, where each is a portion of a plane, i.e., a smooth surface.
General shapes -- in our novel terminology -- arise from $H^1$-deformations of such a prior set $\mathcal{X}_0$. These $H^1$-deformations, evaluated at a prior shape $\Gamma_0=\partial \mathcal{X}_0$, give deformed shapes $\Gi$ if the deformations are injective and continuous. We call these shapes of class $H^{1/2}$ 
and define the set
\begin{equation}
\label{force-ball}
\begin{split}
{\cal H}^{1/2}(\Gamma_0,\R^d):=
\{w\colon \Gamma_0\to X \colon & \exists W\in H^1(\Omega,\Omega) \text{ s.t.}\\ & W\hspace{-.5mm}\,\rule[-2mm]{.1mm}{4mm}_{\hspace{.5mm}\Gamma_0} \text{ injective, continuous}\text{, }W\hspace{-.5mm}\,\rule[-2mm]{.1mm}{4mm}_{\hspace{.5mm}\Gamma_0}=w\}.
\end{split}
\end{equation}
However, in order to have a unique representation for each shape, we have to factor out the homeomorphisms from the prior shape $\Gamma_0$ into itself which are compatible with the set (\ref{force-ball}). Thus, we characterize the following shape space:

\begin{definition}
\label{definition_B12}
Let $\Omega$, $\mathcal{X}_0$ and $\Gamma_0$ be as above. The space of all $H^{1/2}$-shapes is given by
\begin{equation}\label{shape_mainifold}
{\cal B}^{1/2}(\Gamma_0,\R^d):=
{\cal H}^{1/2}(\Gamma_0,\R^d)\big\slash \mbox{\emph{Homeo}}^{1/2}(\Gamma_0),
\end{equation}
where ${\cal H}^{1/2}(\Gamma_0,\R^d)$ is given in (\ref{force-ball}) and $\mbox{\emph{Homeo}}^{1/2}(\Gamma_0)$ is defined by
\begin{equation}
\mbox{\emph{Homeo}}^{1/2}(\Gamma_0)\\ :=
\{ w\colon w\in {\cal H}^{1/2}(\Gamma_0,\R^d)\text{, }w\colon \Gamma_0\to \Gamma_0 \text{ homeomorphism}\}.
\end{equation}
\end{definition}

\begin{remark}
Of course, the properties of the shape space ${\cal B}^{1/2}\left(\Gamma_0,\R^d\right)$ have to be investigated. For example the independence of the prior shape $\Gamma_0$ in the shape space definition is an open question. If it is independent, we can choose, for example, the unit sphere $S^{d-1}$ as prior shape. Another important question is whether the shape space has a manifold structure. Note that this question is very hard and a lot of effort has to be put into it to find the answer. From a theoretical point of view there are several other open questions.
However, this goes beyond the scope of this work and is a topic of subsequent work.
\end{remark}

\begin{remark}
In the following, we assume that ${\cal B}^{1/2}\left(\Gamma_0,\R^d\right)$ has a manifold structure. If necessary, we can refine the space ${\cal B}^{1/2}(\Gamma_0,\R^d)$, e.g., by restriction to an explicit deformation field $W$. In our setting, it arises from the linear elasticity equation and the request of the existence of an arbitrary one is perhaps too strong. This way, we can replace ${\cal H}^{1/2}(\Gamma_0,\R^d)$ by a 
linear space, which is in particular a manifold. However, this conceivable limitation leaves the following theory untouched.
\end{remark}

If $\Gamma\in{\cal B}^{1/2}(\Gamma_0,\R^d)$ is smooth enough to admit a normal vector field $n$, the following isomorphisms naturally arise out of definition (\ref{force-ball}):
\begin{equation}
\label{isomophism}
\begin{split}
& T_\Gamma {\cal B}^{1/2}\left(\Gamma_0,\R^d\right)
\\ & \cong
\{h\colon  h=\phi n  \text{ a.e.}\text{, } \phi\in H^{1/2}(\Gamma) \text{ injective, continuous} \}
\\ & \cong
\{\phi\colon \phi\in H^{1/2}(\Gamma) \text{ injective, continuous} \}
\end{split}
\end{equation}

Now, we can formulate the shape quasi-Newton methods of section \ref{NLP_vol} on the novel shape space ${\cal B}^{1/2}\left(\Gamma_0,\R^d\right)$.
Due to the isomorphism (\ref{isomophism}) and the handy expression (\ref{HadamardConcisely}) of the shape derivative we can state the connection of ${\cal B}^{1/2}\left(\Gamma_0,\R^d\right)$ with respect to $g^S$ to shape calculus, i.e., we can determine a representation $h\in T_{\Gi} {\cal B}^{1/2}(\Gamma_0,\R^d) \cong \{h\colon h\in H^{1/2}(\Gamma) \text{ injective, continuous} \}$ of the shape gradient in terms of $g^S$ defined in (\ref{scp}) by
\begin{equation}
\label{connection1}
g^S(\phi,h)=(r,\phi)_{L^2(\Gi)}
\end{equation}
for all injective and continuous $\phi\in H^{1/2}(\Gi)$,
which is equivalent to
\begin{equation}
\label{connection2}
\int_{\Gi} \phi(s)\cdot [(S^p)^{-1}h](s) \ ds=\int_{\Gi} r(s)\phi(s) \ ds
\end{equation}
for all injective and continuous $\phi\in H^{1/2}(\Gi)$.

Based on the connection (\ref{connection1}) we can formulate the quasi-Newton methods of section \ref{NLP_vol} also on $\left(\mathcal{B}^{1/2}\left(\Gamma_0,\R^d\right),g^S\right)$. From (\ref{connection2}) we get $h=S^pr=(\gamma_0 U)^T n$ where $U\in H^1_0(\Omega,\R^d)$ solves
\begin{equation}
a(U,V)=\int_{\Gi} r\cdot (\gamma_0 V)^T n \ ds
=DJ_{\Gi}[V]=DJ_\Omega[V]\, , \quad \forall\  V\in H^1_0(\Omega,\R^d).
\end{equation}
In general, $h=S^pf=(\gamma_0 U)^\top n$ is not necessarily an element of $T_{\Gi}\mathcal{B}^{1/2}\left(\Gamma_0,\R^d\right)$ because it is not ensured that $U\in H^1_0(\Omega,\mathbb{R}^d)$ is injective and continuous. Under special assumptions depending on the coefficients of a second-order partial differential operator, the right hand-side of a PDE, the domain $\Omega$ on which the PDE is defined and the dimension of $\Omega$, the continuity of a weak solution $U\in H^1_0(\Omega,\mathbb{R}^d)$ of a PDE is guaranteed by the theorem of higher interior regularity \cite[theorem 5, section 6.3]{Evan} combined with the Sobolev embedding theorem. In particular, if these conditions are fulfilled, we get in our two- and three-dimensional case a bounded $C^2$ regularity of $U$. Now, if we require 
\begin{equation}
\label{necessary_cond}
\lVert U\rVert_{\mathcal{C}^1_{b}(\Omega,\mathbb{R}^d)}<1
\end{equation}
we get the injectivity of $U$ due to \cite[lemma 6.13 and remark 6.14]{Allaire}.

\begin{remark}
In implementations, the necessary condition $\lVert U\rVert_{\mathcal{C}^1_{b}(\Omega,\mathbb{R}^d)}<1$ for injectivity of the deformation $U$ is ensured by the particular choice of the Lam\'{e} parameters.
\end{remark}

\end{rev}

\section{Conclusions}
This paper develops an intrinsic metric in shape spaces, which enables to jointly work with domain based and boundary based shape derivative expressions, and which leads to shape optimization algorithms with several computational and analytic advantages as outlined above. 
Furthermore, the metric leads to a novel shape space $\mathcal{B}^{1/2}$.
The properties of $\mathcal{B}^{1/2}$ are beyond the scope of this work and will be touched in subsequent papers.
It is obvious that the results of this paper are not restricted to two space dimensions and also not to interior interface shapes. The whole discussion carries over to shapes which are just parts of the exterior boundary of a computational domain.

\section*{Acknowledgment}
This work has been partly supported by the Deutsche Forschungsgemeinschaft within the Priority program SPP 1648 ``Software for Exascale Computing'' under contract number
Schu804/12-1,  and by BMBF (German Federal Ministry of Education and Research) within the collaborative project R\OE{}NOBIO under contract number 05M13UTA. Furthermore, the authors are indebted to Roland Herzog for many helpful comments and to Vanja Nikoli\'{c} (Klagenfurt University) for triggering the discussion among the authors concerning the novel metric.


\begin{thebibliography}{10}

\bibitem{Absil-book-2008}
{P.-A.} Absil, R.~Mahony, and R.~Sepulchre.
\newblock {\em Optimization algorithms on matrix manifolds}.
\newblock Princeton University Press, 2008.

\bibitem{Agoshkov-1985}
{V. I.} Agoshkov and {V. I.} Lebedev.
\newblock Poincar\'e-Steklov operators and domain decomposition methods in
  variational problems.
\newblock In {\em Computer Procss and Systems}, volume~2, pages 173--227.
  Nauka, Moscow, 1985.
\newblock (in Russian).

\bibitem{Allaire}
G. Allaire.
\newblock {\em Conception optimale de structures}.
\newblock Math\'{e}matiques and Applications 58, Springer, 2007.


\bibitem{Zapletal_2015}
K.~Bandara, F.~Cirak, G.~Of, O.~Steinbach, J.~Zapletal.
\newblock Boundary element based multiresolution shape optimisation in electrostatics.
\newblock {\em Journal of Computational Physics}. 297, 584–598, 2015

\bibitem{Berggren}
M.~Berggren.
\newblock A unified discrete--continuous sensitivity analysis method for shape
  optimization.
\newblock In W.~Fitzgibbon et~al., editor, {\em Applied and numerical partial
  differential equations}, volume~15 of {\em Computational {M}ethods in
  {A}pplied {S}iences}, pages 25--39. Springer, 2010.


\bibitem{CorreaSeger}
R.~Correa and A.~Seeger.
\newblock Directional derivative of a minmax function.
\newblock {\em Nonlinear Analysis}, 9(1):13--22, 1985.

\bibitem{Delfour-Zolesio-2001}
M.~C. Delfour and J.-P. Zol\'esio.
\newblock {\em Shapes and geometries: Analysis, differential calculus, and
  optimization}.
\newblock Advances in Design and Control. SIAM Philadelphia, 2001.

\bibitem{EpplerHarbrechtSchneider-2007}
K. Eppler, H. Harbrecht and R. Schneider.
\newblock  On convergence in elliptic shape optimization.
\newblock {\em SIAM Journal on Control and Optimization.} 46(1):61--83, 2007.

\bibitem{Evan}
L. C. Evans.
\newblock {\em Partial Differential Equations}.
\newblock Amer. Math. Soc., Providence, RI, 1993.

\bibitem{Langer-2015}
P. Gangl, U. Langer, A. Laurain, H. Meftahi, and K. Sturm.
\newblock Shape optimization of an electric motor subject to nonlinear
  magnetostatics.
\newblock Technical report, \url{http://arxiv.org/abs/1501.04752}, 2015.

\bibitem{GrossReusken}
S.~Gross and A.~Reusken.
\newblock {\em Numerical methods for two-phase incompressible flows}, volume~40
  of {\em Computational Mathematics}.
\newblock Springer, 2010.

\bibitem{HipPag_2015}
R.~Hiptmair and A. Paganini.
\newblock Shape Optimization by Pursuing
Diffeomorphisms.
\newblock {\em Comput. Methods Appl. Math.} 15(3):291–305, 2015



\bibitem{KhW-2004}
B. N. Khoromskij and G. Wittum.
\newblock {\em Numerical solution of elliptic differential equations by
  reduction to the interface}.
\newblock Number~36 in Lecture Notes in Computational Science and Engineering.
  Springer, 2004.
  
\bibitem{LaurainSturm2013}
A.~Laurain and K. Sturm.
\newblock Domain expression of the shape derivative and application to electrical impedance tomography.
\newblock Technical Report No. 1863, WIAS Berlin, 2013.


\bibitem{meyer2003discrete}
M.~Meyer, M.~Desbrun, P.~Schr{\"o}der, and A.~H. Barr.
\newblock Discrete differential-geometry operators for triangulated
  2-manifolds.
\newblock In {\em Visualization and mathematics III}, pages 35--57. Springer,
  2003.

\bibitem{MM-2006}
{P. W.} Michor and D.~Mumford.
\newblock Riemannian geometries on spaces of plane curves.
\newblock {\em Journal of the European Mathematical Society}, 8:1--48, 2006.

\bibitem{neagel20015scalable}
A.~N\"agel, V. H.~Schulz, M.~Siebenborn, and G.~Wittum.
\newblock Scalable methods for structured inverse modelling in diffusive
  processes.
\newblock {\em Computing and Visualization in Science}, 2015 (submitted).

\bibitem{Skin-2015}
A. N\"agel, V. H. Schulz, M. Siebenborn, and G. Wittum.
\newblock Scalable shape optimization methods for structured inverse modeling
  in 3d diffusive processes.
\newblock {\em Computing and Visualization in Science}, 2015.

\bibitem{Novruzi-2002}
A. Novruzi and M. Pierre.
\newblock Structure of shape derivatives.
\newblock {\em Journal of Evolution Equations}, 2:365--382, 2002.

\bibitem{Paganini}
A.~Paganini.
\newblock Approximative shape gradients for interface problems.
\newblock Technical Report 2014-12, Seminar for Applied Mathematics, ETH
  Z{\"u}rich, 2014.

\bibitem{SS-2009}
S.~Schmidt and V. H.~Schulz.
\newblock Impulse response approximations of discrete shape {H}essians with
  application in {C}{F}{D}.
\newblock {\em SIAM Journal on Control and Optimization}, 48(4):2562--2580,
  2009.

\bibitem{AIAA-2013}
S. Schmidt, C. Ilic, V. H. Schulz, and N. Gauger.
\newblock Three-dimensional large-scale aerodynamic shape optimization based on
  shape calculus.
\newblock {\em {AIAA} Journal}, 51(11):2615--2627, 2013.

\bibitem{VHS-shape-Riemann}
{V. H.} Schulz.
\newblock A {R}iemannian view on shape optimization.
\newblock {\em Foundations of Computational Mathematics}, 14:483--501, 2014.

\bibitem{Schulz-Structure-2014}
V. H. Schulz, M. Siebenborn, and K. Welker.
\newblock Structured inverse modeling in parabolic diffusion problems.
\newblock {\em {SIAM} Journal on Control and Optimization}, 53(6): 3319--3338, 2015.
\newblock \url{http://arxiv.org/abs/1409.3464}.

\bibitem{Schulz-LN-2014}
V. H. Schulz, M. Siebenborn, and K. Welker.
\newblock Towards a {L}agrange-{N}ewton approach for {P}{D}{E} constrained
  shape optimization.
\newblock In A. Pratelli and G. Leugering, editors, {\em New Trends in Shape Optimization}, volume 166 of {\em
  International Series of Numerical Mathematics}, pp 229--249. Springer, 2015.
\newblock \url{http://arxiv.org/abs/1405.3266}.

\bibitem{SokoZol}
J.~Sokolowski and {J.-P.} Zol\'{e}sio.
\newblock {\em An introduction to shape optimization}.
\newblock Springer, 1992.

\bibitem{Bletzinger-2015}
E. Stavropoulou, M. Hojjat, and K.-U. Bletzinger.
\newblock In-plane mesh regularization for node-based shape optimization
  problems.
\newblock {\em Computer Methods in Applied Mechanics and Engineering}, 275:39--54, 2014.

\bibitem{Sturm2013}
K. Sturm.
\newblock Lagrange method in shape optimization for non-linear partial
  differential equations: A material derivative free approach.
\newblock Technical Report No. 1817, WIAS Berlin, 2013.

\bibitem{Troeltzsch}
F.~Tr\"oltzsch.
\newblock {\em {Optimal control of partial differential equations: Theory,
  methods, and applications}}, volume 112 of {\em Applied Mathematics}.
\newblock American Mathematical Society, 2010.

\bibitem{Berggren-horn-2007}
R.~Udawalpola and M.~Berggren.
\newblock Optimization of an acoustic horn with respect to efficiency and
  directivity.
\newblock {\em International Journal for Numerical Methods in Engineering}, 73(11):1571--1606,
  2007.

\bibitem{Sokolowski-1996}
J. Sokolowski.
\newblock Displacement Derivatives in Shape Optimization of Thin Shells.
\newblock Research Report RR-2995, INRIA, pp.21, 1996.


\end{thebibliography}
\bibliographystyle{plain}

\end{document}